\documentclass[a4paper,11pt,draft]{amsart}
\usepackage{amssymb}

\newtheorem{theorem}{Theorem}
\newtheorem{thm}[theorem]{Theorem}

\theoremstyle{definition}

\theoremstyle{remark}
\newtheorem{remark}[theorem]{Remark}

\newcommand{\BB}{{\mathbb B}}
\newcommand{\DD}{{\mathbb D}}

\newcommand{\EE}{{\mathcal E}}
\newcommand{\CC}{{\mathbb C}}

\renewcommand{\phi}{\varphi}

\begin{document}

\title{Graphs of multifunctions}

\author{Armen Edigarian}

\address{Institute of Mathematics, Jagiellonian University,
Reymonta 4/526, 30-059 Krak\'ow, Poland}
\email{Armen.Edigarian@im.uj.edu.pl}
\thanks{The author was supported in part by the KBN grant
No. 5 P03A 033 21}

\date{26 January 2004}


\maketitle


\section{Introduction}

In \cite{Shch1} N.~Shcherbina gives a positive answer to the
following question of Nishino (see \cite{Nishino1}):

\textit{Let $\DD$ be the unit disc in $\CC$ and let $f:\DD\to\CC$
be a continuous function such that its graph
\begin{equation}
\Gamma(f)=\{\big(z,f(z)\big): z\in\DD\}
\end{equation}
is a pluripolar subset of $\CC^2$. Does ti follow that $f$ is
holomorphic?}

Recently, in \cite{Shch2} N.~Shcherbina proves a more general
result for multi-values functions defined by Weierstrass
pseudopolynomial.

\begin{thm}\label{thm:1} Let $D$ be a domain in $\CC^n$ and let
$\EE$ be a subset of $D\times\CC$ of the form
\begin{equation}\label{eq:1}
\EE=\{(z,w)\in D\times\CC: w^m+a_1(z)w^{m-1}+\dots+a_m(z)=0\},
\end{equation}
where $a_1,\dots,a_m$ are continuous complex-valued functions on
$D$. Then $\EE$ is a pluripolar subset of $\CC^{n+1}$ if and only
if the functions $a_1,\dots,a_m$ are holomorphic on $D$.
\end{thm}

The proof of Theorem~\ref{thm:1}, presented in \cite{Shch2}, is
based on Oka's result (see e.g.~\cite{Nishino2}, Theorem~4.9) and
is a generalization of the methods of \cite{Shch1}.

The main purpose of the paper is to show that Theorem~\ref{thm:1}
is a simple corollary of the results of \cite{Shch1}. Actually,
our proof is similiar, in spirit, to the proof of the Oka's
result.

\section{Proofs}

Recall the following result of Rad\'o (see e.g.~\cite{Nar}).
\begin{thm} Let $D$ be a domain in $\CC^n$ and let $f$ be a
continuous function on $D$. Assume that $f$ is holomorphic on
$D\setminus f^{-1}\{0\}$. Then $f$ is holomorphic on $D$.
\end{thm}

\begin{proof}[Proof of Theorem~\ref{thm:1}]
It suffices to show that from pluripolarity of $\EE$ it follows
holomorphicity of $a_1,\dots,a_m$. We prove by induction on $m$.

For $m=1$ it's the main result in \cite{Shch1}.

So, assume that Theorem~\ref{thm:1} holds for any domain $D$ and
for $1,2,\dots,m-1$. Put
\begin{equation*}
P(z,w)=w^m+a_1(z)w^{m-1}+\dots+a_m(z).
\end{equation*}
Fix $z_0\in D$. Let $w_1,\dots,w_s$ be all distinct solutions of
the equation $P(z_0,\cdot)=0$. Assume that $s>1$ (i.e. the
multiplicity at any point $w_j$ is less than $m$).

Let us show that $a_1,\dots,a_m$ are holomorphic in a neighborhood
of $z_0$. Indeed, fix $j\in\{1,\dots,s\}$. Let $k$ be the
multiplicity of $P(z_0,\cdot)=0$ at $w_j$ (we know that $k<m$).
Then there exists an $r>0$ such that for any
$z\in\BB(z_0,r)=\{\zeta\in\CC^n:\|\zeta-z_0\|<r\}$ in
$\DD(w_j,r)=\{\xi\in\CC: |\xi-w_j|<r\}$ there is exactly $k$
solutions of the equation $P(z,\cdot)=0$ (counted with
multiplicity). Moreover, there is no solution on
$\partial\DD(w_j,r)$. Assume that $\zeta_1(z),\dots,\zeta_k(z)$
are these solutions. Put
\begin{equation*}
\widetilde P(z,w)=(w-\zeta_1(z))\dots(w-\zeta_k(z)).
\end{equation*}
Note $\widetilde P$ is well-defined on $\BB(z_0,r)$ and that its
coefficients are continuous in $\BB(z_0,r)$. Since
\begin{equation*}
\{(z,w)\in\BB(z_0,r)\times\CC: \widetilde P(z,w)=0\} \subset \EE,
\end{equation*}
by induction step we get that its coefficients are holomorphic in
$\BB(z_0,r)$.

We may repeat the same argument for all the solutions
$w_1,\dots,w_s$ of the equation $P(z_0,\cdot)=0$ and get that
$a_1,\dots,a_m$ (as polynomials of the previously obtained
coefficients) are holomorphic in some neighborhood of $z_0$.

Let $B$ denote the set of all points $z\in D$ such that
$P(z,\cdot)=0$ has a solution with the multiplicity $m$. Then $B$
is a closed subset of $D$. Moreover, $a_1,\dots,a_m$ are
holomorphic on $D\setminus B$. If $B$ is a pluripolar set, then
$a_1,\dots,a_m$ are holomorphic on $D$ (use the continuity). So,
assume that $B$ is non-pluripolar.

Let
\begin{equation*}
\phi_k(z,w)=\frac{\partial^{k}P(z,w)}{\partial w^{k}}(z,w).
\end{equation*}
Note that
$\phi_{m-2}(z,w)=\frac{m!}{2}w^2+(m-1)!a_1(z)w+(m-2)!a_2(z)$. Let
\begin{equation*}
\Delta(z)=\big((m-1)!a_1(z)\big)^2-4\frac{m!}{2}(m-2)!a_2(z).
\end{equation*}
Note that $B\subset\{z\in D: \Delta(z)=0\}$. Moreover, $\Delta$ is
a continuous function on $D$ and is holomorphic on $D\setminus
\Delta^{-1}\{0\}$. Hence, by Rad\'o's theorem, $\Delta$ is a
holomorphic function on $D$. Since $B$ is non-pluripolar, we have
$\Delta\equiv0$. So,
$a_2(z)=\frac{m(m-1)}{2}\Big(\frac{a_1(z)}{2}\Big)^2$.

Let us show by (inverse) induction on $\ell$, that
\begin{equation*}
\phi_\ell(z,w)=\frac{m!}{(m-\ell)!}\Big(w-\frac{a_1(z)}{m}\Big)^{m-\ell}
\end{equation*}
for any $\ell=1,\dots,m$. We know that it's true for $\ell=m, m-1,
m-2$. Assume that it holds for $m,m-1,\dots,\ell+1$. We want to
show that it holds for $\ell$. By induction step
\begin{multline*}
\phi_\ell(z,w)=\frac{m!}{(m-\ell)!}\Big(w-\frac{a_1(z)}{m}\Big)^{m-\ell}+
(m-\ell)!a_\ell(z)\\
+(-1)^{m-\ell+1}\frac{m!}{(m-\ell)!}\Big(\frac{a_1(z)}{m}\Big)^{m-\ell}.
\end{multline*}
Put $\psi_\ell=
(m-\ell)!a_\ell(z)+(-1)^{m-\ell+1}\frac{m!}{(m-\ell)!}\Big(\frac{a_1(z)}{m}\Big)^{m-\ell}$.
Note that $\psi_\ell$ is a continuous function and $B\subset\{z\in
D:\psi_\ell(z)=0\}$. Moreover, it is holomorphic on
$D\setminus\psi_\ell^{-1}\{0\}$. Hence, $\psi_\ell$ is holomorphic
on $D$ (again use Rado's theorem) and, therefore,
$\psi_\ell\equiv0$.

So, $B=D$ and $P(z,w)=(w-\frac{a_1(z)}{n})^m$. By induction step,
$a_1$ (hence, $a_2,\dots,a_m$) is a holomorphic function.
\end{proof}

\begin{remark} Note that if in Theorem~\ref{thm:1} there exists a
point $z_0\in D$ such that $P(z_0,\cdot)$ has $m$ distinct
solutions, then one can give even simpler proof. Indeed, put
$\Delta(z)=\prod_{j\ne k}\big(\xi_j(z)-\xi_k(z)\big)$, where
$\xi_1(z),\dots,\xi_{m}(z)$ are the solutions (counted with
multiplicity) of $P(z,\cdot)=0$. Then $\Delta$ is a symmetric
polynomial of $\xi_1(z),\dots,\xi_m(z)$, so it is a polynomial of
$a_1(z),\dots,a_m(z)$. Therefore, $\Delta$ is a continuous
function on $D$. Note that $\Delta$ and $a_1,\dots,a_m$ are
holomorphic functions on $D\setminus\Delta^{-1}\{0\}$ (use the
results of \cite{Shch1}). We know that $\Delta(z_0)\ne0$. So,
$\Delta^{-1}(\{0\})$ is a proper analytic set. Since
$a_1,\dots,a_m$ are continuous on $D$ and holomorphic on
$D\setminus\Delta^{-1}\{0\}$, they are holomorphic on~$D$.
\end{remark}

\begin{remark} N.~Shcherbina informed me that he also found,
independently, similar proof.
\end{remark}

\bibliographystyle{amsplain}


\end{document}